\documentclass[reqno,12pt,a4paper]{amsart}

\usepackage{mathrsfs}
\usepackage{amssymb}
\usepackage{amsfonts}
\usepackage{latexsym}
\usepackage{amsthm}
\usepackage{graphicx}

\def\lb{\label}

\newcommand{\er}[1]{\textrm{(\ref{#1})}}

\begin{document}

%%%%%%%%%% Some definitions %%%%%%%%%%

%%%%%%%% Equations, theorems %%%%%%%%%
\renewcommand{\theequation}{\arabic{section}.\arabic{equation}}
\theoremstyle{plain}
\newtheorem{theorem}{\bf Theorem}[section]
\newtheorem{lemma}[theorem]{\bf Lemma}
\newtheorem{corollary}[theorem]{\bf Corollary}
\newtheorem{proposition}[theorem]{\bf Proposition}
\newtheorem{definition}[theorem]{\bf Definition}
\newtheorem{example}[theorem]{\bf Example}

\theoremstyle{remark}
\newtheorem*{remarks}{\bf Remarks}
\newtheorem*{remark}{\bf Remark}
%\newtheorem{example}{\bf Example}

%%%%% Alphabet %%%%%
\def\a{\alpha}  \def\cA{{\mathcal A}}     \def\bA{{\bf A}}  \def\mA{{\mathscr A}}
\def\b{\beta}   \def\cB{{\mathcal B}}     \def\bB{{\bf B}}  \def\mB{{\mathscr B}}
\def\g{\gamma}  \def\cC{{\mathcal C}}     \def\bC{{\bf C}}  \def\mC{{\mathscr C}}
\def\G{\Gamma}  \def\cD{{\mathcal D}}     \def\bD{{\bf D}}  \def\mD{{\mathscr D}}
\def\d{\delta}  \def\cE{{\mathcal E}}     \def\bE{{\bf E}}  \def\mE{{\mathscr E}}
\def\D{\Delta}  \def\cF{{\mathcal F}}     \def\bF{{\bf F}}  \def\mF{{\mathscr F}}
\def\c{\chi}    \def\cG{{\mathcal G}}     \def\bG{{\bf G}}  \def\mG{{\mathscr G}}
\def\z{\zeta}   \def\cH{{\mathcal H}}     \def\bH{{\bf H}}  \def\mH{{\mathscr H}}
\def\e{\eta}    \def\cI{{\mathcal I}}     \def\bI{{\bf I}}  \def\mI{{\mathscr I}}
\def\p{\psi}    \def\cJ{{\mathcal J}}     \def\bJ{{\bf J}}  \def\mJ{{\mathscr J}}
\def\vT{\Theta} \def\cK{{\mathcal K}}     \def\bK{{\bf K}}  \def\mK{{\mathscr K}}
\def\k{\kappa}  \def\cL{{\mathcal L}}     \def\bL{{\bf L}}  \def\mL{{\mathscr L}}
\def\l{\lambda} \def\cM{{\mathcal M}}     \def\bM{{\bf M}}  \def\mM{{\mathscr M}}
\def\L{\Lambda} \def\cN{{\mathcal N}}     \def\bN{{\bf N}}  \def\mN{{\mathscr N}}
\def\m{\mu}     \def\cO{{\mathcal O}}     \def\bO{{\bf O}}  \def\mO{{\mathscr O}}
\def\n{\nu}     \def\cP{{\mathcal P}}     \def\bP{{\bf P}}  \def\mP{{\mathscr P}}
\def\r{\rho}    \def\cQ{{\mathcal Q}}     \def\bQ{{\bf Q}}  \def\mQ{{\mathscr Q}}
\def\s{\sigma}  \def\cR{{\mathcal R}}     \def\bR{{\bf R}}  \def\mR{{\mathscr R}}
\def\S{\Sigma}  \def\cS{{\mathcal S}}     \def\bS{{\bf S}}  \def\mS{{\mathscr S}}
\def\t{\tau}    \def\cT{{\mathcal T}}     \def\bT{{\bf T}}  \def\mT{{\mathscr T}}
\def\f{\phi}    \def\cU{{\mathcal U}}     \def\bU{{\bf U}}  \def\mU{{\mathscr U}}
\def\F{\Phi}    \def\cV{{\mathcal V}}     \def\bV{{\bf V}}  \def\mV{{\mathscr V}}
\def\P{\Psi}    \def\cW{{\mathcal W}}     \def\bW{{\bf W}}  \def\mW{{\mathscr W}}
\def\o{\omega}  \def\cX{{\mathcal X}}     \def\bX{{\bf X}}  \def\mX{{\mathscr X}}
\def\x{\xi}     \def\cY{{\mathcal Y}}     \def\bY{{\bf Y}}  \def\mY{{\mathscr Y}}
\def\X{\Xi}     \def\cZ{{\mathcal Z}}     \def\bZ{{\bf Z}}  \def\mZ{{\mathscr Z}}

\def\vr{\varrho}
%*********************
\def\be{{\bf e}} \def\bc{{\bf c}}
\def\bv{{\bf v}} \def\bu{{\bf u}}
\def\Om{\Omega} \def\bp{{\bf p}}
%************************
\def\bbD{\pmb \Delta}
\def\mm{\mathrm m}
\def\mn{\mathrm n}
%*************************

\newcommand{\mc}{\mathscr {c}}

\newcommand{\gA}{\mathfrak{A}}          \newcommand{\ga}{\mathfrak{a}}
\newcommand{\gB}{\mathfrak{B}}          \newcommand{\gb}{\mathfrak{b}}
\newcommand{\gC}{\mathfrak{C}}          \newcommand{\gc}{\mathfrak{c}}
\newcommand{\gD}{\mathfrak{D}}          \newcommand{\gd}{\mathfrak{d}}
\newcommand{\gE}{\mathfrak{E}}
\newcommand{\gF}{\mathfrak{F}}           \newcommand{\gf}{\mathfrak{f}}
\newcommand{\gG}{\mathfrak{G}}           %\newcommand{\gg}{\mathfrak{g}}
\newcommand{\gH}{\mathfrak{H}}           \newcommand{\gh}{\mathfrak{h}}
\newcommand{\gI}{\mathfrak{I}}           \newcommand{\gi}{\mathfrak{i}}
\newcommand{\gJ}{\mathfrak{J}}           \newcommand{\gj}{\mathfrak{j}}
\newcommand{\gK}{\mathfrak{K}}            \newcommand{\gk}{\mathfrak{k}}
\newcommand{\gL}{\mathfrak{L}}            \newcommand{\gl}{\mathfrak{l}}
\newcommand{\gM}{\mathfrak{M}}            \newcommand{\gm}{\mathfrak{m}}
\newcommand{\gN}{\mathfrak{N}}            \newcommand{\gn}{\mathfrak{n}}
\newcommand{\gO}{\mathfrak{O}}
\newcommand{\gP}{\mathfrak{P}}             \newcommand{\gp}{\mathfrak{p}}
\newcommand{\gQ}{\mathfrak{Q}}             \newcommand{\gq}{\mathfrak{q}}
\newcommand{\gR}{\mathfrak{R}}             \newcommand{\gr}{\mathfrak{r}}
\newcommand{\gS}{\mathfrak{S}}              \newcommand{\gs}{\mathfrak{s}}
\newcommand{\gT}{\mathfrak{T}}             \newcommand{\gt}{\mathfrak{t}}
\newcommand{\gU}{\mathfrak{U}}             \newcommand{\gu}{\mathfrak{u}}
\newcommand{\gV}{\mathfrak{V}}             \newcommand{\gv}{\mathfrak{v}}
\newcommand{\gW}{\mathfrak{W}}             \newcommand{\gw}{\mathfrak{w}}
\newcommand{\gX}{\mathfrak{X}}               \newcommand{\gx}{\mathfrak{x}}
\newcommand{\gY}{\mathfrak{Y}}              \newcommand{\gy}{\mathfrak{y}}
\newcommand{\gZ}{\mathfrak{Z}}             \newcommand{\gz}{\mathfrak{z}}

\def\ve{\varepsilon}   \def\vt{\vartheta}    \def\vp{\varphi}    \def\vk{\varkappa}

\def\A{{\mathbb A}} \def\B{{\mathbb B}} \def\C{{\mathbb C}}
\def\dD{{\mathbb D}} \def\E{{\mathbb E}} \def\dF{{\mathbb F}} \def\dG{{\mathbb G}} \def\H{{\mathbb H}}\def\I{{\mathbb I}} \def\J{{\mathbb J}} \def\K{{\mathbb K}} \def\dL{{\mathbb L}}\def\M{{\mathbb M}} \def\N{{\mathbb N}} \def\O{{\mathbb O}} \def\dP{{\mathbb P}} \def\R{{\mathbb R}}\def\S{{\mathbb S}} \def\T{{\mathbb T}} \def\U{{\mathbb U}} \def\V{{\mathbb V}}\def\W{{\mathbb W}} \def\X{{\mathbb X}} \def\Y{{\mathbb Y}} \def\Z{{\mathbb Z}}

%%%%% Arrows %%%%%

\def\la{\leftarrow}              \def\ra{\rightarrow}            \def\Ra{\Rightarrow}
\def\ua{\uparrow}                \def\da{\downarrow}
\def\lra{\leftrightarrow}        \def\Lra{\Leftrightarrow}

%%%%% Typography %%%%%

\def\lt{\biggl}                  \def\rt{\biggr}
\def\ol{\overline}               \def\wt{\widetilde}
\def\ul{\underline}
\def\no{\noindent}

%%%%% Math signs %%%%%

\let\ge\geqslant                 \let\le\leqslant
\def\lan{\langle}                \def\ran{\rangle}
\def\/{\over}                    \def\iy{\infty}
\def\sm{\setminus}               \def\es{\emptyset}
\def\ss{\subset}                 \def\ts{\times}
\def\pa{\partial}                \def\os{\oplus}
\def\om{\ominus}                 \def\ev{\equiv}
\def\iint{\int\!\!\!\int}        \def\iintt{\mathop{\int\!\!\int\!\!\dots\!\!\int}\limits}
\def\el2{\ell^{\,2}}             \def\1{1\!\!1}
\def\sh{\sharp}
\def\wh{\widehat}
\def\bs{\backslash}
\def\intl{\int\limits}
%%%%% Math operations %%%%%

\def\na{\mathop{\mathrm{\nabla}}\nolimits}
\def\sh{\mathop{\mathrm{sh}}\nolimits}
\def\ch{\mathop{\mathrm{ch}}\nolimits}
\def\where{\mathop{\mathrm{where}}\nolimits}
\def\all{\mathop{\mathrm{all}}\nolimits}
\def\as{\mathop{\mathrm{as}}\nolimits}
\def\Area{\mathop{\mathrm{Area}}\nolimits}
\def\arg{\mathop{\mathrm{arg}}\nolimits}
\def\const{\mathop{\mathrm{const}}\nolimits}
\def\det{\mathop{\mathrm{det}}\nolimits}
\def\diag{\mathop{\mathrm{diag}}\nolimits}
\def\diam{\mathop{\mathrm{diam}}\nolimits}
\def\dim{\mathop{\mathrm{dim}}\nolimits}
\def\dist{\mathop{\mathrm{dist}}\nolimits}
\def\Im{\mathop{\mathrm{Im}}\nolimits}
\def\Iso{\mathop{\mathrm{Iso}}\nolimits}
\def\Ker{\mathop{\mathrm{Ker}}\nolimits}
\def\Lip{\mathop{\mathrm{Lip}}\nolimits}
\def\rank{\mathop{\mathrm{rank}}\limits}
\def\Ran{\mathop{\mathrm{Ran}}\nolimits}
\def\Re{\mathop{\mathrm{Re}}\nolimits}
\def\Res{\mathop{\mathrm{Res}}\nolimits}
\def\res{\mathop{\mathrm{res}}\limits}
\def\sign{\mathop{\mathrm{sign}}\nolimits}
\def\span{\mathop{\mathrm{span}}\nolimits}
\def\supp{\mathop{\mathrm{supp}}\nolimits}
\def\Tr{\mathop{\mathrm{Tr}}\nolimits}
\def\BBox{\hspace{1mm}\vrule height6pt width5.5pt depth0pt \hspace{6pt}}

%%%%%%%%%%%%% specialities %%%%%%%%%%%%%%

\newcommand\nh[2]{\widehat{#1}\vphantom{#1}^{(#2)}}
%{{\mathop{#1}\limits^\wedge}\vphantom{#1}^{(#2)}}
\def\dia{\diamond}

\def\Oplus{\bigoplus\nolimits}

%%%%%%%%%%% End of definitions %%%%%%%%%%

%%%%% OLD OLD OLD

\def\qqq{\qquad}
\def\qq{\quad}
\let\ge\geqslant
\let\le\leqslant
\let\geq\geqslant
\let\leq\leqslant
\newcommand{\ca}{\begin{cases}}
\newcommand{\ac}{\end{cases}}
\newcommand{\ma}{\begin{pmatrix}}
\newcommand{\am}{\end{pmatrix}}
\renewcommand{\[}{\begin{equation}}
\renewcommand{\]}{\end{equation}}
\def\eq{\begin{equation}}
\def\qe{\end{equation}}
\def\[{\begin{equation}}
\def\bu{\bullet}

\title[On spectrum of magnetic Schr\"odinger operators on periodic graphs]
{On continuous spectrum of magnetic Schr\"odinger operators on periodic discrete graphs}

\date{\today}
\author[Evgeny Korotyaev]{Evgeny Korotyaev}
\address{Saint-Petersburg State University, Universitetskaya nab. 7/9, St. Petersburg, 199034, Russia,
\ korotyaev@gmail.com, \
e.korotyaev@spbu.ru,}
\author[Natalia Saburova]{Natalia Saburova}
\address{Northern (Arctic) Federal University, Severnaya Dvina Emb. 17, Arkhangelsk, 163002, Russia,
 \ n.saburova@gmail.com, \ n.saburova@narfu.ru}


\begin{abstract}
We consider Schr\"odinger operators with periodic electric and
magnetic potentials on periodic discrete graphs. The spectrum of
such operators consists of an absolutely continuous (a.c.) part (a~union of
a finite number of non-degenerate bands) and a finite number of eigenvalues  of infinite multiplicity. We prove the following results:
1) the a.c. spectrum of the magnetic Schr\"odinger operators is empty for
specific graphs and magnetic fields; 2) we obtain necessary and
sufficient conditions under which the a.c. spectrum of the magnetic
Schr\"odinger operators is empty; 3) the spectrum of the magnetic Schr\"odinger operator with each
magnetic potential $t\a$, where $t$ is a coupling constant, has an
a.c. component for all except finitely many $t$ from any bounded
interval.
\end{abstract}

\maketitle

\section {\lb{Sec1}Introduction and main results}
\setcounter{equation}{0}

There are many physical phenomena described by periodic
Schr\"odinger operators. In general, the spectral analysis of a
periodic operator $H$ is based on the so-called Floquet
decomposition, where this operator has a representation  as a direct
integral of a family of fiber operators $H(k)$ with discrete
spectra, see e.g., \cite{RS78}. The fiber operator $H(k)$ depends on
the so-called quasimomentum $k$ belonging to the torus
$\T^d=\R^d/(2\pi \Z)^d$. Its eigenvalues $\l_j(k)$,
$j=1,2,\ldots$\,, arranged in non-decreasing order, depend on $k$
continuously. The spectrum of the operator $H$ is a union of
spectral bands $\s_j$, arising as the ranges of the \emph{band
functions} $\l_j(\cdot)$, i.e., $\s_j=\l_j(\T^d)$. If some function
$\l_j(k)=\L_j=\const$ on an open domain of $\T^d$, then $\L_j$ is an
eigenvalue of $H$ of infinite multiplicity. We call $\{\L_j\}$ a
\emph{flat band} or a \emph{degenerate} band. All other bands are
\emph{non-degenerate}. The union of all flat bands is the \emph{flat
spectrum} of $H$ and the union of all non-degenerate bands is the
a.c. spectrum of $H$. One of important problems is to describe
the spectrum: when the operator $H$ has the flat spectrum and when
it does not. There are a lot of results about the spectrum of
periodic differential operators, see, e.g., \cite{BS98,D85,S99,T73}.

There are some results about the spectrum of discrete Schr\"odinger
operators  with periodic electric potentials on periodic graphs, see
\cite{HN09,KS14,KS19}. It is shown in \cite{KS19} that the first
band of these operators is non-degenerate. Thus, the spectrum of the
Schr\"odinger operators always has an a.c. component. Some classes of periodic graphs on which the discrete Laplacian has only the a.c. spectrum are described in \cite{HN09}.

We consider discrete Schr\"odinger operators with periodic magnetic
and electric potentials on periodic graphs. This case and even the
case of magnetic Laplacians are much more complicated compared to
the non-magnetic one. We do not know results about the flat spectrum
of the discrete magnetic Schr\"odinger operators. Our goal is to
show that, in contrast to the non-magnetic case, the a.c. spectrum
of the magnetic Schr\"odinger operators may be empty (for specific
graphs and magnetic potentials) and to determine necessary and
sufficient conditions for absence of the a.c. spectrum. Moreover, we
show that absence of the a.c. spectrum is a quite rare situation,
i.e., the spectrum of the magnetic Schr\"odinger operator  with a
periodic magnetic potential $t\a$, where $t\in \R$ is a coupling
constant, has an a.c. component for all except finitely many $t$ from
any bounded interval.

\subsection{Magnetic Schr\"odinger operators on periodic graphs}
Let $\cG=(\cV,\cE)$ be a connected infinite graph, possibly  having
loops and multiple edges and embedded into the space $\R^d$. Here
$\cV$ is  the set of its vertices and $\cE$ is the set of its
unoriented edges. Considering each edge in $\cE$ to have two
orientations, we introduce the set $\cA$ of all oriented edges. An
edge starting at a vertex $u$ and ending at a vertex $v$ from $\cV$
will be denoted as the ordered pair $(u,v)\in\cA$ and is said to be
\emph{incident} to the vertices. Let $\ul\be=(v,u)$ be the inverse
edge of $\be=(u,v)\in\cA$. We define the degree ${\vk}_v$ of the vertex $v\in\cV$ as the number of
all edges from $\cA$, starting at $v$.

Let $\G$ be a lattice of rank $d$ in $\R^d$ with a basis $\{a_1,\ldots,a_d\}$, i.e.,
$$
\G=\Big\{a : a=\sum_{s=1}^dn_sa_s, \; n_s\in\Z,\; s\in\N_d\Big\}, \qqq \N_d=\{1,\ldots,d\},
$$
and let
\[\lb{fuce}
\Omega=\Big\{x\in\R^d : x=\sum_{s=1}^dt_sa_s, \; 0\leq t_s<1,\; s\in\N_d\Big\}
\]
be the \emph{fundamental cell} of the lattice $\G$. We define the equivalence relation on $\R^d$:
$$
x\equiv y \; (\hspace{-4mm}\mod \G) \qq\Leftrightarrow\qq x-y\in\G \qqq
\forall\, x,y\in\R^d.
$$

We consider \emph{locally finite $\G$-periodic graphs} $\cG$, i.e., graphs satisfying the
following conditions: \\
$\bu$\; $\cG=\cG+a$ for any $a\in\G$ and  the quotient graph
$\cG_*=\cG/\G$ is finite.

The basis $a_1,\ldots,a_d$ of the lattice $\G$ is called the {\it
periods} of $\cG$.  We call the quotient graph $\cG_*=\cG/\G$ the
\emph{fundamental graph} of the periodic graph $\cG$. The
fundamental graph $\cG_*$ is a graph on the $d$-dimensional torus
$\R^d/\G$. The graph $\cG_*=(\cV_*,\cE_*)$ has the vertex set
$\cV_*=\cV/\G$, the set $\cE_*=\cE/\G$ of unoriented edges and the
set $\cA_*=\cA/\G$ of oriented edges which are finite.

Let $\ell^2(\cV)$ be the Hilbert space of all square summable
functions $f:\cV\to \C$ equipped with the norm
$$
\|f\|^2_{\ell^2(\cV)}=\sum_{v\in\cV}|f(v)|^2<\infty.
$$
We consider the magnetic Schr\"odinger operators $H_\a$ acting on $\ell^2(\cV)$ and given by
\[\lb{SOH1}
H_\a=\D_\a+Q,
\]
where $Q$ is an electric potential and   $\D_\a$ is the magnetic
Laplacian having the form
\begin{equation}
\lb{col}
(\D_\a f)(v)=\sum\limits_{\be=(v,\,u)\in\cA}\big(f(v)-e^{i\a(\be)}f(u)\big),
\qqq f\in \ell^2(\cV), \qqq v\in \cV,
\end{equation}
$\a:\cA\ra\R$ is a \emph{magnetic vector potential} on $\cG$, i.e., it satisfies
\[\lb{MPot}
\a(\ul\be\,)=-\a(\be),\qqq \forall\,\be\in\cA.
\]
The sum in \er{col} is taken over all edges from $\cA$ starting at
the vertex $v$.  We assume that the magnetic potential $\a$ and the
electric potential $Q$ are $\G$-periodic, i.e., they satisfy
\[\lb{Gpps}
\a(\be+a)=\a(\be),\qqq Q(v+a)=Q(v),\qqq \forall\,(v,\be,a)\in\cV\ts\cA\ts\G.
\]
It is known  that the magnetic Schr\"odinger operator $H_\a$ is a
bounded self-adjoint operator on $\ell^2(\cV)$ (see, e.g.,
\cite{HS99}).

\subsection{Spectrum of the magnetic Schr\"odinger operator.}
The magnetic Schr\"odinger operator $H_\a$ on $\ell^2(\cV)$ has the
decomposition into a constant fiber direct integral given by
\[
\lb{raz}
\mH={1\/(2\pi)^d}\int^\oplus_{\T^d}\ell^2(\cV_*)\,dk ,\qqq
UH_\a U^{-1}={1\/(2\pi)^d}\int^\oplus_{\T^d}H_\a(k)dk,
\]
$\T^d=\R^d/(2\pi\Z)^d$, for some unitary operator
$U:\ell^2(\cV)\to\mH$.  Here $\ell^2(\cV_*)=\C^\nu$ is the fiber
space, $\nu=\#\cV_*$. The parameter $k\in\T^d$ is called the
\emph{quasimomentum}. The precise expression of the fiber operator
$H_\a(k)$ is given by \er{Hvt'} -- \er{fado}. Note that $H_\a(0)$ is
the magnetic Schr\"odinger operator on $\cG_*$. Each fiber operator
$H_\a(k)$, $k\in\T^d$, has $\n$ real eigenvalues $\l_{\a,j}(k)$,
$j\in\N_\n$, which are labeled in non-decreasing order
(counting multiplicities) by
\[\label{eq.3}
\l_{\a,1}(k)\leq\l_{\a,2}(k)\leq\ldots\leq\l_{\a,\nu}(k),
\qqq \forall\,k\in\T^d.
\]
Each $\l_{\a,j}(\cdot)$, $j\in\N_\n$, is a real and piecewise
analytic function  on the torus $\T^d$ and creates the
\emph{spectral band} $\s_j(H_\a)$ given by
\[
\lb{ban.1}
\s_{j}(H_\a)=[\l_{\a,j}^-,\l_{\a,j}^+]=\l_{\a,j}(\T^d).
\]

Note that if $\l_{\a,j}(\cdot)=\L_j=\const$ on some subset of $\T^d$
of positive Lebesgue measure, then  the operator $H_\a$ on $\cG$ has
the eigenvalue $\L_{j}$ of infinite multiplicity. We call $\{\L_j\}$
a \emph{flat band}. Thus, the spectrum of the operator $H_\a$ on the
periodic graph $\cG$ is given by
\[\lb{spec}
\s(H_\a)=\bigcup_{k\in\T^d}\s\big(H_\a(k)\big)=
\bigcup_{j=1}^{\nu}\s_j(H_\a)=\s_{ac}(H_\a)\cup \s_{fb}(H_\a),
\]
where $\s_{ac}(H_\a)$ is the a.c. spectrum, which
is a union of non-degenerate bands, and $\s_{fb}(H_\a)$ is the \emph{flat
spectrum} which is the set of all flat bands (eigenvalues of
infinite multiplicity).

It is known (see, e.g., Proposition 4.2 in \cite{HN09}) that $\l_*$
is  an eigenvalue of infinite multiplicity of the operator $H_\a$ if
and only if $\l_*$ is an eigenvalue of the operator $H_\a(k)$ for
all $k\in\T^d$. It gives another labeling of the eigenvalues of
$H_\a(k)$. If the operator $H_\a$ has $r\ge 0$ eigenvalues of
infinite multiplicity
\[\lb{fb2}
\m_1\le \m_2\le\ldots\le \m_r, \qqq r\le \n,
\]
then we can separate constant eigenvalues $\m_j(\cdot)=\m_j=\const$,
of $H_\a(\cdot)$, $j\in\N_r$. All other eigenvalues $\m_j(k)$,
$r<j\le \n$, of $H_\a(k)$ can be enumerated in non-decreasing order
 by
\[\label{eq.3.1}
\m_{r+1}(k)\leq\m_{r+2}(k)\leq\ldots\leq\m_{\n}(k), \qqq
\forall\,k\in\T^d,
\]
counting multiplicities. We define the \emph{spectral bands}
$\gS_j(H_\a)=[\m_j^-,\m_j^+]=\m_j(\T^d)$, $j\in\N_{\n}$, where each
band $\gS_j(H_\a)$, $r<j\le \n$, is non-degenerate, i.e.,
$\m_j^-<\m_j^+$. Thus, the number of non-degenerate spectral bands
of the operator $H_\a$ is $\n-r$.

\subsection{Main results}
Recall that the a.c. spectrum of the Schr\"odinger operator $H_0$ with
the magnetic potential $\a=0$ is not empty, see \cite{KS19}.
The following simple examples show that this is not true for the
magnetic Schr\"odinger operator $H_\a$, i.e., the a.c. spectrum of
$H_\a$ on $\cG$ is empty for specific periodic graphs and magnetic
potentials $\a$. In the first example we construct the magnetic
Schr\"odinger operators $H_\a$ when their fiber operators $H_\a(k)$
do not depend on the quasimomentum $k\in\T^d$. It is clear that in
this case the spectrum of $H_\a$ is flat, since all band functions
are constants.

\begin{example}\lb{Ex1}
Let each oriented edge of a periodic graph $\cG$ have multiplicity
2, and let a periodic magnetic vector potential $\a$ on $\cG$
satisfy
\[\lb{mpco1}
|\a(\be)-\a(\wt\be\,)|=\pi\qq \textrm{for each pair of multiple edges}\;
(\be,\wt\be\,)\in\cA^2.
\]
Then for any periodic electric potential $Q$ the magnetic Schr\"odinger operator $H_\a=\D_\a+Q$ has  the
decomposition \er{raz}, where its fiber operator has the form
$H_\a(k)=\diag\big(\vk_v+Q(v)\big)_{v\in\cV_*}$, and $\vk_v$ is the degree
of the vertex $v\in\cV_*$.
In particular, the spectrum of $H_\a$ is flat and is given by
$\big\{\vk_v+Q(v)\big\}_{v\in\cV_*}$.
\end{example}

The second  example shows that the condition $H_\a(\cdot)=\const$ is
not necessary for absence of the a.c. spectrum
of $H_\a$.

\setlength{\unitlength}{1.0mm}
\begin{figure}[h]
\centering
\unitlength 1.2mm % = 2.845pt
\linethickness{0.4pt}
\ifx\plotpoint\undefined\newsavebox{\plotpoint}\fi % GNUPLOT compatibility
\begin{picture}(115,23)(0,0)
\put(0,10){\line(1,0){5.00}}
\put(20,10){\line(1,0){15.00}}
\put(5,10){\circle*{1}}
\put(20,10){\circle*{1}}
\bezier{200}(5,10)(12.5,0)(20,10)
\bezier{200}(5,10)(12.5,20)(20,10)
\put(10,15.5){$\scriptstyle\be_{-2,2}$}
\put(10,3){$\scriptstyle\be_{-2,1}$}
\put(26,11){$\scriptstyle\be_{-1,1}$}
\put(41,15.5){$\scriptstyle\be_{0,2}$}
\put(41,3){$\scriptstyle\be_{0,1}$}
\put(71,15.5){$\scriptstyle\be_{2,2}$}
\put(71,3){$\scriptstyle\be_{2,1}$}
\put(56,11){$\scriptstyle\be_{1,1}$}
\put(50,10){\line(1,0){15.00}}
\put(35,10){\circle*{1}}
\put(50,10){\circle*{1}}
\bezier{200}(35,10)(42.5,0)(50,10)
\bezier{200}(35,10)(42.5,20)(50,10)
\put(65,10){\circle*{1}}
\put(80,10){\circle*{1}}
\bezier{200}(65,10)(72.5,0)(80,10)
\bezier{200}(65,10)(72.5,20)(80,10)
\put(80,10){\line(1,0){5.00}}
\put(19.5,7.5){$\scriptstyle-1$}
\put(2,7.5){$\scriptstyle-2$}
\put(50,7.5){$\scriptstyle1$}
\put(34,7.5){$\scriptstyle0$}
\put(80,7.5){$\scriptstyle3$}
\put(64,7.5){$\scriptstyle2$}
\put(0,15){$\cG$}

\put(0,4){\emph{a})}
%*********************************
\put(98,15){$\cG_*$}
\put(96,4){\emph{b})}
\put(100,10){\circle*{1}}
\put(115,10){\circle*{1}}
\bezier{200}(100,10)(107.5,0)(115,10)
\bezier{200}(100,10)(107.5,20)(115,10)
\put(115,8){$\scriptstyle v_1$}
\put(98.5,8){$\scriptstyle v_0$}
\put(100,10){\line(1,0){15.00}}
\put(107,15){\vector(1,0){1.00}}
\put(107,10){\vector(1,0){1.00}}
\put(107,5){\vector(1,0){1.00}}
\put(107,8.0){$\scriptstyle \be_3$}
\put(107,16){$\scriptstyle\be_2$}
\put(107,3){$\scriptstyle\be_1$}

\end{picture}

\caption{\footnotesize  \emph{a}) The periodic graph $\cG=(\Z,\cE)$; \; \emph{b}) the fundamental graph $\cG_*=\cG/(2\Z)$.} \label{Fig1}
\end{figure}

\begin{example}\lb{Ex2}
We consider the periodic graph $\cG=(\Z,\cE)$, where the edge set is given by
$$
\cE=\big\{\be_{n,1}=(n,n+1)\; \textrm{for all} \; n\in\Z
\big\}\cup\big\{\be_{n,2}=(n,n+1)\; \textrm{for all even} \; n\in\Z\big\}
$$
(see Fig.\ref{Fig1}a). Let $H_\a=\D_\a+Q$ be the Schr\"odinger
operator  on $\cG$, where a 2-periodic electric potential $Q$ and a
periodic magnetic potential $\a$ are given by
\[\lb{gamm}
\ca Q(0)=0 \\
Q(1)\in \R \ac,\qq
\a(\be_{n,s})=\ca \a_o,  & \textrm{if $s=1$ and $n$ is even}\\
\pi+\a_o,  & \textrm{if $s=2$ and $n$ is even}\\
0, & \textrm{if $s=1$ and $n$ is odd}\ac,\qq n\in\Z, \qq s=1,2,
\]
and $\a_o\in \R$. Then the fiber operator $H_\a(k)$ depends on
$k\in\T^d$, but the spectrum of $H_\a$ is flat and is given by
\[\lb{fsE1}
\textstyle\l_{s}=3+\frac12\,\big(Q(1)+(-1)^s\sqrt{Q^2(1)+4}\,\big),
\qqq s=1,2.
\]
\end{example}

The proof of Examples \ref{Ex1} and \ref{Ex2} is given in Subsection
\ref{SSecE}.

\begin{remark}
Below in Theorem \ref{Tfbs} we determine necessary and sufficient
conditions when the spectrum of the magnetic Schr\"odinger operator $H_\a$
is flat, i.e., the a.c. spectrum of $H_\a$ is empty. In order to present this result we need more
definitions.
\end{remark}

Examples \ref{Ex1} and \ref{Ex2} show that all spectral bands of the
magnetic Schr\"odinger  operator $H_\a$ may be
flat. In the following theorem we prove that
this is a quite rare situation.

\begin{theorem}\lb{Tfmp}
Let $H_{t\a}=\D_{t\a}+Q$ be the magnetic Schr\"odinger operator
defined  by \er{SOH1} -- \er{col} with a periodic magnetic potential
$t\a$ and a periodic electric potential $Q$ on a periodic graph
$\cG$, where $t$ is a coupling constant. Then  the
a.c. spectrum of $H_{t\a}$ is not empty for all except
finitely many $t\in[0,1]$.
\end{theorem}

\begin{remarks}
1) It is known \cite{BS98,S99} that the spectrum of the
Schr\"odinger operator with periodic electric and magnetic
potentials on $L^2(\R^d)$ is a.c. and the proof is based on the
Thomas' approach.
We do not use this argument. In order to prove Theorem \ref{Tfmp} we
determine trace formulas for the magnetic Schr\"odinger operator
(see Theorem \ref{TFNL0}) as Fourier series and use classical
function theory. Moreover, as an unperturbed case we use the result
from \cite{KS19} that the first spectral band of the Schr\"odinger
operator without magnetic field is not flat.

2) Theorem \ref{Tfmp} also holds true for weighted magnetic Laplace
and Schr\"odinger operators, in particular, for the normalized ones.
The proof repeats the proof of Theorem \ref{Tfmp}.
\end{remarks}

\subsection{Historical review}
There are a lot papers about the spectrum of periodic differential
operators,  see articles \cite{BS98,D85,HH95,S99,T73}  and the
references therein. The first result in this direction was obtained
by Thomas \cite{T73}. He proved that each band function of
Schr\"odinger operators with periodic potentials on $L^2(\R^d)$ is
not flat. Later on this approach was used in many papers
\cite{BS98,D85,HH95,S99}.
Danilov \cite{D85} proved that the spectrum of the Dirac operator with
periodic potentials on $\R^d$ is a.c. Hempel and
Herbst \cite{HH95} proved that the spectrum of magnetic Laplacians
with small periodic magnetic vector potentials  in $\R^d$ is a.c.
Birman and Suslina \cite{BS98} (for the case $d=2$), and Sobolev
\cite{S99} (for $d\ge 2$) proved that the spectrum of the
Schr\"odinger operator with periodic electric and magnetic
potentials on $L^2(\R^d)$ is a.c.

In the discrete settings the situation is quite different. The
spectrum of the discrete Schr\"odinger operator with periodic
electric and magnetic potentials on periodic graphs consists of a
finite number of bands. Some of them may be degenerate. Thus,  the
spectrum of the discrete Schr\"odinger operators has an a.c.
component which is a union of all non-degenerate bands and a flat
component which is the set of all degenerate bands (eigenvalues of
infinite multiplicity). In \cite{KS19} it was proved that the first
band of the discrete Schr\"odinger operators with periodic electric
potentials on periodic graphs is always non-degenerate, i.e., the
a.c. spectrum is not empty. It was also shown  that all except two
bands of the Schr\"odinger operators may be degenerate (see
Proposition 7.2 in \cite{KS14}) and the number of the flat bands
 may be arbitrary. We do not know results about the flat spectrum
of the discrete magnetic Schr\"odinger operators.

\section{Proof}\lb{Sec2}
\setcounter{equation}{0}

\subsection{Edge indices}
For each $x\in\R^d$ we introduce the vector $x_\A\in\R^d$ by
\[\lb{cola}
x_\A=(x_1,\ldots,x_d), \qqq \textrm{where} \qq x=
\textstyle\sum\limits_{s=1}^dx_sa_s,
\]
i.e., $x_\A$ is the coordinate vector of $x$ with respect  to the
basis $\A=\{a_1,\ldots,a_d\}$ of the lattice $\G$.

For any vertex $v\in\cV$ of a $\G$-periodic graph $\cG$ the
following unique representation holds true:
\[\lb{Dv}
v=v_0+[v], \qq \textrm{where}\qq v_0\in\Omega,\qquad [v]\in\G,
\]
$\Omega$ is a fundamental cell of the lattice $\G$ defined by
\er{fuce}. In other words, each vertex $v$ can be obtained from a
vertex  $v_0\in \Omega$ by a shift by a vector $[v]\in\G$. For any
oriented edge $\be=(u,v)\in\cA$ we define the \emph{edge index}
$\t(\be)$ as the vector of the lattice $\Z^d$ given by
\[
\lb{in}
\t(\be)=[v]_\A-[u]_\A\in\Z^d,
\]
where $[v]\in\G$ is defined by \er{Dv} and the vector $[v]_\A\in\Z^d$ is given by \er{cola}.

The edge indices $\t(\be)$ and the magnetic potential $\a$ are $\G$-periodic, i.e., they satisfy
\[\lb{Gpe}
\t(\be+a)=\t(\be),\qqq \a(\be+a)=\a(\be),\qqq \forall\, (\be,a)\in\cA \ts\G.
\]
On the set $\cA$ of all oriented edges of the $\G$-periodic graph $\cG$ we define the surjection
\[\lb{sur}
\gf:\cA\rightarrow\cA_*=\cA/\G,
\]
which maps each $\be\in\cA$ to its equivalence class
$\be_*=\gf(\be)$  which is an oriented edge of the fundamental graph
$\cG_*$. The identities \er{Gpe} and the mapping $\gf$ allow us to
define uniquely

$\bu$ the \emph{magnetic potential} $\a$ on edges of the fundamental
graph $\cG_*=(\cV_*,\cA_*)$ which is induced by the magnetic
potential $\a$:
\[
\lb{Mdco}
\a(\be_*)=\a(\be)\qq \textrm{ for some $\be\in\cA$ \; such that }  \;
 \be_*=\gf(\be), \qqq \be_*\in\cA_*;
\]

$\bu$ \emph{indices} of the fundamental graph edges which are induced by the
indices of the periodic graph edges:
\[
\lb{dco}
\t(\be_*)=\t(\be) \qq \textrm{ for some $\be\in\cA$ \; such that }  \;
 \be_*=\gf(\be), \qqq \be_*\in\cA_*.
\]

\subsection{Cycle indices and magnetic fluxes} A \emph{path} $\bp$
 in a graph $\cG=(\cV,\cA)$ is a sequence of consecutive edges
\[\lb{depa}
\bp=(\be_1,\be_2,\ldots,\be_n),
\]
where $\be_s=(v_{s-1},v_s)\in\cA$, $s=1,\ldots,n$, for some vertices
$v_0,v_1,\ldots,v_n\in\cV$. The vertices $v_0$ and $v_n$ are called
the \emph{initial} and \emph{terminal} vertices of the path $\bp$,
respectively. If $v_0=v_n$, then the path $\bp$ is called a
\emph{cycle}. The number $n$ of edges in a cycle $\bc$ is called
the \emph{length} of $\bc$ and is denoted by $|\bc|$, i.e.,
$|\bc|=n$.

\begin{remark} A path $\bp$ is uniquely defined by the sequence of it's oriented
edges $(\be_1,\be_2,\ldots,\be_n)$. The sequence of it's vertices
$(v_0,v_1,\ldots,v_n)$ does not uniquely define~$\bp$, since
multiple edges are allowed in the graph $\cG$.
\end{remark}

Let $\cC$ be the sets of all cycles of the fundamental graph $\cG_*$.
For any cycle $\bc\in\cC$ we define

$\bu$ the \emph{cycle index} $\t(\bc)\in\Z^d$ by
\[\lb{cyin}
\t(\bc)=\sum\limits_{\be\in\bc}\t(\be),  \qqq \bc\in\cC,
\]

$\bu$  the \emph{flux} $\a(\bc)\in[-\pi,\pi]$ of the magnetic potential
 $\a$ through the cycle $\bc$ by
\[\lb{cyin1}
\a(\bc)=\bigg(\sum\limits_{\be\in\bc}\a(\be)\bigg)\hspace{-3mm} \mod 2\pi,
 \qqq \qqq \bc\in\cC.
\]
Note that we consider fluxes of
the magnetic potential $\a$ modulo $2\pi$, since the magnetic
potential $\a$ appears in the Laplacian $\D_\a$ via the factor
$e^{i\a(\be)}$, $\be\in\cA$. For each cycle $\bc=(\be_1,\ldots,\be_n)$ we define  a cycle
$\ul\bc=(\ul\be_n,\ldots,\ul\be_1)$. From the definition of indices and
fluxes it follows that
\[\lb{inin}
\begin{array}{l}
\t(\ul\be)=-\t(\be), \qqq \a(\ul\be)=-\a(\be), \qqq \forall\,\be\in\cA_*; \\[6pt]
\t(\ul\bc)=-\t(\bc), \qqq \a(\ul\bc)=-\a(\bc), \qqq  \,\forall\,\bc\in\cC.
\end{array}
\]

\subsection{Direct integral} We introduce the Hilbert space, i.e., a constant fiber direct
integral,
\[\lb{Hisp}
\mH=L^2\Big(\T^{d},{dk\/(2\pi)^d}\,,\ell^2(\cV_*)\Big)=\int_{\T^{d}}^{\os}\ell^2(\cV_*)\,{dk
\/(2\pi)^d}\,,
%\qqq \cH=\ell^2(\cV_*),
\qqq \T^d=\R^d/(2\pi\Z)^d,
\]
 equipped with the norm $
\|g\|^2_{\mH}=\int_{\T^d}\|g(k,\cdot)\|_{\ell^2(\cV_*)}^2\frac{dk}{(2\pi)^d}\,,
$ where the function $g(k,\cdot)\in\ell^2(\cV_*)$ for almost all
$k\in\T^d$. We identify the vertices of the fundamental
graph $\cG_*=(\cV_*,\cE_*)$ with the vertices  from the fundamental
cell $\Omega$. The unitary Gelfand transform $U:\ell^2(\cV)\to\mH$
is given by
\[
\lb{UGt}
(Uf)(k,v)=\sum\limits_{\mm=(m_1,\ldots,m_d)\in\Z^d}e^{-i\lan
\gm,k\ran } f(v+m_1a_1+\ldots+m_da_d), \qqq (k,v)\in \T^d\ts\cV_*,
\]
where $\{a_1,\ldots,a_d\}$ is the basis of the lattice $\G$, and
$\lan\cdot\,,\cdot\,\ran$ denotes the standard inner product in
$\R^d$. We recall Theorem 3.1 from \cite{KS17}.

\begin{theorem}\label{TFDC}
The magnetic Schr\"odinger operator $H_\a=\D_\a+Q$ on $\ell^2(\cV)$
has the following decomposition into a constant fiber direct
integral
\[
\lb{raz1}
\begin{aligned}
& UH_\a U^{-1}=\int^\oplus_{\T^d}H_\a(k){dk\/(2\pi)^d}\,,
\end{aligned}
\]
where $U:\ell^2(\cV)\to\mH$ is the unitary Gelfand transform defined
by \er{UGt}, and the fiber magnetic Schr\"odinger operator $H_\a(k)$
on $\ell^2(\cV_*)$ has the form
\[\label{Hvt'}
H_\a(k)=\D_\a(k)+Q, \qqq \forall\,k\in\T^d.
\]
Here $Q$ is the electric potential on $\ell^2(\cV_*)$ and $\D_\a(k)$
is the fiber magnetic Laplacian given by
\[
\label{fado}
\big(\D_\a(k)f\big)(v)=\sum_{\be=(v,\,u)\in\cA_*} \big(f(v)
-e^{i(\a(\be)+\lan\t(\be),\,k\ran)}f(u)\big), \qq f\in\ell^2(\cV_*),\qq v\in \cV_*,
\]
where $\t(\be)$ is the index of the edge $\be\in\cA_*$ defined by
\er{in}, \er{dco}.
\end{theorem}

From Theorem \ref{TFDC} it follows that  the fiber magnetic
Schr\"odinger operator $H_\a(k)$
 is  the $\n\ts\n$ matrix  defined by \er{Hvt'} -- \er{fado} and  in the standard orthonormal
basis of $\ell^2(\cV_*)=\C^\n$, $\n=\#\cV_*$, is given by
\[\lb{mrSo}
H_\a(k)=q-A_\a(k), \qqq q=\diag(q_v)_{v\in\cV_*}, \qqq q_v=\vk_v+Q(v),
\]
where $\vk_v$ is the degree of the vertex $v$, and the matrix
$A_\a(k)=\big(A_{\a,uv}(k)\big)_{u,v\in\cV_*}$ has the form
\[\lb{mrmL}
A_{\a,uv}(k)=\sum\limits_{\be=(u,v)\in\cA_*}e^{-i(\a(\be)+\lan\t(\be),\,k\ran)}.
\]

\begin{lemma}\lb{Lfbs}
Let $H_\a(k)$ be the fiber operator given by \er{mrSo}, \er{mrmL}. Then $\l_*$ is an eigenvalue of infinite multiplicity of the magnetic Schr\"odinger operator $H_\a$ if and only if $\l_*$ is an eigenvalue of $H_\a(k)$ for all $k\in\T^d$.
\end{lemma}

\no\textbf{Proof.} The proof is quite standard. But for the reader's convenience we repeat it. Let $\l_*$ be an eigenvalue of infinite multiplicity of $H_\a$. Then $\det\big(\l_* I_\n-H_\a(k)\big)=0$ for all $k\in\cB$, where $\cB$ is a subset of $\T^d$ of positive Lebesgue measure, and $I_\n$ is the identity $\n\ts\n$ matrix. The function $f(k)=\det\big(\l_* I_\n-H_\a(k)\big)$ is real analytic in $k\in\T^d$ (moreover, it is an entire function of $k\in\C^{d}$). Then $f(k)=0$ for all $k\in\T^d$, i.e., $\l_*$ is an eigenvalue of $H_\a(k)$ for all $k\in\T^d$. The converse is obvious. \qq \BBox

\subsection{Proof of Examples.} \lb{SSecE}
We prove Examples \ref{Ex1} and \ref{Ex2} where we construct
magnetic Schr\"odinger operators with the empty a.c. spectrum.

\medskip

\no \textbf{Proof of Example \ref{Ex1}.} Let each oriented edge of a periodic graph $\cG$ have multiplicity 2,  and let a magnetic vector potential $\a$ on $\cG$ satisfy \er{mpco1}. Then for each edge $\be=(u,v)\in\cA_*$ of the fundamental graph $\cG_*=(\cV_*,\cA_*)$ there exists an edge $\wt\be=(u,v)\in\cA_*$ such that
$$
\t(\wt\be\,)=\t(\be),\qq\textrm{and}\qq |\a(\be)-\a(\wt\be\,)|=\pi.
$$
Thus, the matrix
$A_\a(k)=\big(A_{\a,uv}(k)\big)_{u,v\in\cV_*}$ given by \er{mrmL} has the form
\begin{multline*}
\textstyle A_{\a,uv}(k)=\sum\limits_{\be=(u,v)\in\cA_*}e^{-i(\a(\be)+\lan\t(\be),\,k\ran)}
=\frac12\sum\limits_{\be=(u,v)\in\cA_*}\big(e^{-i(\a(\be)+\lan\t(\be),\,k\ran)}+
e^{-i(\a(\wt\be\,)+\lan\t(\wt\be\,),\,k\ran)}\big)\\\textstyle=
\frac12\sum\limits_{\be=(u,v)\in\cA_*}e^{-i\lan\t(\be),\,k\ran}
\big(e^{-i\a(\be)}+e^{-i(\a(\be)\pm\pi)}\big)=0.
\end{multline*}
This and \er{mrSo} yield that the fiber magnetic Schr\"odinger
operator $H_\a(k)$ is diagonal and has the form
$
H_\a(k)=\diag\big(\vk_v+Q(v)\big)_{v\in\cV_*}.
$
Thus, the spectrum of
$H_\a$  on $\cG$ is flat and is given by
$\big\{\vk_v+Q(v)\big\}_{v\in\cV_*}$. \qq \BBox

\medskip

\no \textbf{Proof of Example \ref{Ex2}.} The fundamental graph $\cG_*=\cG/(2\Z)$ consists of two vertices $v_0,v_1$ and three multiple edges $\be_1,\be_2,\be_3$ connecting these vertices (Fig.~\ref{Fig1}\emph{b}) with indices
$$
\t(\be_1)=\t(\be_2)=0,\qqq \t(\be_3)=1,
$$
and their inverse edges. The magnetic potential $\a$ on edges of $\cG_*$  is given by
$$
\a(\be_1)=\a_o,\qqq \a(\be_2)=\pi+\a_o,\qqq \a(\be_3)=0,
$$
for some $\a_o\in \R$. Then the fiber magnetic Schr\"odinger
operator $H_{\a}(k)$, $k\in\T$,  given by \er{mrSo} -- \er{mrmL}, on
$\cG_*$ has the form
$$
H_{\a}(k)=\ma
   3 & -e^{-i\a_0}-e^{-i(\pi+\a_o)}-e^{-ik} \\
  -e^{i\a_0}-e^{i(\pi+\a_0)}-e^{ik} & 3+Q(1)
\am =\ma    3 & -e^{-ik} \\   -e^{ik} & 3+Q(1) \am,
$$
where $Q(0)=0$. Thus, the eigenvalues of $H_{\a}(k)$  are given by
\er{fsE1}. \qq \BBox

\subsection{Trace formulas} In order to formulate trace formulas for
the fiber magnetic Schr\"odinger operator $H_\a(k)$, we need some
modifications of the fundamental graph $\cG_*$. We add a loop
$\be_v$ with index $\t(\be_v)=0$ and the magnetic potential $\a(\be_v)=0$ at each vertex $v$ of the
fundamental graph $\cG_*=(\cV_*,\cA_*)$ and consider the modified
fundamental graph $\wt\cG_*=(\cV_*,\wt\cA_*)$,  where
\[\lb{wtAs}
\wt\cA_*=\cA_*\cup\{\be_v\}_{v\in\cV_*}.
\]
We denote by $\wt\cC$ the set of all cycles on $\wt\cG_*$.
For each cycle $\bc\in\wt\cC$ we define the \emph{weight}
\[\lb{Wcy}
\o(\bc)=\o(\be_1)\ldots \o(\be_n), \qq \textrm{where}\qq \bc=(\be_1,\ldots,\be_n)\in\wt\cC,
\]
and $\o(\be)$ is defined by
\[\lb{webe}
\o(\be)=\left\{
\begin{array}{cl}
-1,  & \qq \textrm{if} \qq  \be\in\cA_* \\[2pt]
q_v, & \qq \textrm{if} \qq \be=\be_v
\end{array}\right.,\qqq q_v=\vk_v+Q(v).
\]

The next lemma shows that the operator $H_\a(k)$ can be considered as a
fiber weighted magnetic operator on the modified fundamental graph $\wt\cG_*=(\cV_*,\wt\cA_*)$.

\begin{lemma}\lb{LDeH}
The fiber magnetic Schr\"odinger operator $H_\a(k)=\big(H_{\a,uv}(k)\big)_{u,v\in\cV_*}$
given by \er{mrSo}, \er{mrmL} satisfies
\[\lb{mrmL+}
H_{\a,uv}(k)=\sum\limits_{\be=(u,v)\in\wt\cA_*}\o(\be) e^{-i(\a(\be)+\lan\t(\be),\,k\ran)},\qqq
\forall\, u,v\in\cV_*, \qqq \forall\,k\in\T^d,
\]
where $\wt\cA_*$ is the set of all edges of the modified fundamental
graph $\wt\cG_*$ defined by \er{wtAs}; $\o(\be)$ is given by
\er{webe}, and $\t(\be)$  is the index of the edge $\be\in\cA_*$
defined by \er{in}, \er{dco}.
\end{lemma}

\no \textbf{Proof.} Let $u,v\in\cV_*$. If $u\neq v$, then, using
\er{webe}  and \er{mrSo}, \er{mrmL}, we have
$$
\sum\limits_{\be=(u,v)\in\wt\cA_*}\o(\be) e^{-i(\a(\be)+\lan\t(\be),\,k\ran)}
=-\sum\limits_{\be=(u,v)\in\cA_*} e^{-i(\a(\be)+\lan\t(\be),\,k\ran)}=H_{\a,uv}(k).
$$
Similarly, if $u=v$, then we obtain
\begin{multline*}
\sum\limits_{\be=(v,v)\in\wt\cA_*}\o(\be) e^{-i(\a(\be)+\lan\t(\be),\,k\ran)}
=\o(\be_v)-\sum\limits_{\be=(v,v)\in\cA_*} e^{-i(\a(\be)+\lan\t(\be),\,k\ran)}\\
=q_v-\sum\limits_{\be=(v,v)\in\cA_*} e^{-i(\a(\be)+\lan\t(\be),\,k\ran)}=H_{\a,vv}(k).
\end{multline*}
Thus, the identity \er{mrmL+} has been proved. \qq \BBox

\medskip

Let $\wt\cC_{n,\gm}$ be the set of all cycles from $\wt\cC$ of
length $n$ and with index $\gm\in\Z^d$:
\[\lb{cNnm}
\wt\cC_{n,\gm}=\{\bc\in\wt\cC: |\bc|=n
\;\textrm{and}\;\t(\bc)=\gm\},
\]
where $|\bc|$ is the length of the cycle $\bc$, and $\t(\bc)$ is
the index of $\bc$ defined by \er{cyin}.

In the following theorem we determine all Fourier coefficients of
$\Tr H_\a^n(k)$ as functions of $k\in \T^d$. This is a crucial
point for our consideration.

\begin{theorem}\lb{TFNL0} Let $H_\a(k)$, $k\in\T^d$, be the fiber
magnetic Schr\"odinger operator defined by \er{Hvt'} -- \er{fado}.
Then for each $n\in\N$ the trace of $H_\a^n(k)$ has the following
Fourier series representation
\[
 \lb{FsTrD}
 \begin{aligned}
& \Tr H_\a^n(k)=
\sum\limits_{\gm\in\Z^d}\gh_{\a,n,\gm}e^{-i\lan\gm,k\ran},
 \\
&  \gh_{\a,n,\gm}=\sum_{\bc\in\wt\cC_{n,\gm}}\o(\bc)e^{-i\a(\bc)},\qq
\supp \gh_{\a,n,\,\cdot}\ss \{\gm\in\Z^d: \|\gm\|\le n\t_+  \},
\end{aligned}
\]
where $\t_+=\max_{\be\in\cA_*}\|\t(\be)\|$, $\t(\be)$  is the
index of the edge $\be\in\cA_*$, and $\|\cdot\|$ is the standard norm in $\R^d$. Here $\wt\cC_{n,\gm}$ is defined by
\er{cNnm}; $\a(\bc)$ is the flux of the magnetic potential $\a$
through the cycle $\bc$ defined by \er{cyin}, and $\o(\bc)$ is given
by \er{Wcy}.
\end{theorem}

\begin{remark}
1) The formulas \er{FsTrD} are \emph{trace formulas}, where the
traces of the fiber operators are expressed in terms of some
geometric parameters of the graph (vertex degrees, cycle indices and
lengths).

2)  There are trace formulas  for Schr\"odinger operators with
periodic potentials on the line, see e.g., \cite{K97}. They were
used to obtain two-sided estimates of potentials in terms of gap
lengths (or a solution for KdV in terms of the action variables) in
\cite{K00} via the conformal mapping theory for the quasimomentum.
Unfortunately, we do not know results about trace formulas for the
multidimensional case.
\end{remark}

\no\textbf{Proof.} Using \er{mrmL+}, for each $n\in\N$ we obtain
$$
\begin{aligned}
& \Tr H_\a^n(k)=
\sum_{v_1,\ldots,v_n\in\cV_*}H_{\a,v_1v_2}(k)H_{\a,v_2v_3}(k)\ldots
 H_{\a,v_{n-1}v_n}(k)H_{\a,v_nv_1}(k)
 \\
 & =
\sum_{v_1,\ldots,v_n\in\cV_*} \sum_{\be_1,\ldots
\be_n\in\wt\cA_*}\o(\be_1)\o(\be_2)\ldots \o(\be_n)
e^{-i(\a(\be_1)+\a(\be_2)+\ldots+\a(\be_n)+
\lan\t(\be_1)+\t(\be_2)+\ldots+\t(\be_n),k\ran)}
\\
& =
\sum_{\bc\in\wt\cC_n}\o(\bc)e^{-i(\a(\bc)+\lan\t(\bc),k\ran)},\qqq
\textrm{where}\qqq  \be_j=(v_j,v_{j+1}), \qq j\in \N_n,\
v_{n+1}=v_1,
\end{aligned}
$$
 and $\wt\cC_n$ is the set of all cycles of length $n$ on
$\wt\cG_*$. Thus, we have  the finite Fourier series for
$2\pi\Z^d$-periodic function $\Tr H_\a^n(k)$, since $\t(\bc)\in
\Z^d$. We rewrite this Fourier series in the standard form
\[\lb{rg11+}
\begin{aligned}
\Tr
H_\a^n(k)=\sum_{\bc\in\wt\cC_n}\o(\bc)e^{-i(\a(\bc)+\lan\t(\bc),k\ran)}
=\sum\limits_{\gm\in\Z^d}\sum_{\bc\in\wt\cC_{n,\gm}}
\o(\bc)e^{-i\a(\bc)}e^{-i\lan\gm,k\ran}
\\
=\sum\limits_{\gm\in\Z^d}e^{-i\lan\gm,k\ran}\sum_{\bc\in\wt\cC_{n,\gm}}\o(\bc)e^{-i\a(\bc)}
=\sum\limits_{\gm\in\Z^d}\gh_{\a,n,\gm}e^{-i\lan\gm,k\ran},
\end{aligned}
\]
where the coefficients $\gh_{\a,n,\gm}$ have the form
\[\lb{ghe}
\gh_{\a,n,\gm}=\sum_{\bc\in\wt\cC_{n,\gm}}\o(\bc)e^{-i\a(\bc)}.
\]

By the definition of the cycle index \er{cyin}, for each cycle $\bc$ of length $n$ we have
$$
\textstyle\|\t(\bc)\|\leq \sum\limits_{\be\in\bc}\|\t(\be)\|\leq n\t_+, \qqq\textrm{where}\qqq \t_+=\max\limits_{\be\in\cA_*}\|\t(\be)\|.
$$
Thus, $\gh_{\a,n,\gm}=0$ for all $\gm\in\Z^d$ such that $\|\gm\|>n\t_+$. \qq
 \BBox

\subsection{Proof of the main Theorems}
We present necessary and sufficient conditions under which the
spectrum of the magnetic Schr\"odinger operators is flat.

\begin{theorem}\lb{Tfbs}
Let $H_\a=\D_\a+Q$ be the magnetic Schr\"odinger operator defined by
\er{SOH1} -- \er{col} with  a periodic magnetic potential $\a$ and
a periodic electric potential $Q$ on a periodic graph
$\cG$. Then the following statements are equivalent:

(i) The spectrum of $H_\a$ is flat.

(ii) For each $n\in\N_\n$ the trace of $H_\a^n(k)$ does not depend on $k\in\T^d$.

(iii) The Fourier coefficients $\gh_{\a,n,\gm}$ of $\Tr H_\a^n(k)$ for
all $(n,\gm)\in\N_\n\ts\big(\Z^d\sm\{0\}\big)$ satisfy
\[\lb{mpco}
\gh_{\a,n,\gm}=\sum_{\bc\in\wt\cC_{n,\gm}}\o(\bc)e^{-i\a(\bc)}=0.
\]

(iv) For any $n\in\N_\n$ the following identities hold true
\[
\lb{Trsq} \frac1{(2\pi)^d}\int_{\T^d}\Tr^2
H_\a^n(k)dk=|\gh_{\a,n,0}|^2,\qq {\where}\qq \gh_{\a,n,0}=
\sum_{\bc\in\wt\cC_{n,0}}\o(\bc)e^{-i\a(\bc)}.
\]
Here $\wt\cC_{n,\gm}$ is defined by \er{cNnm}; $\a(\bc)$ is the flux
of the magnetic potential $\a$ through the cycle $\bc$ defined by
\er{cyin}, and $\o(\bc)$ is given by \er{Wcy}.
\end{theorem}

\begin{remark}
Theorem \ref{Tfbs} also determines the necessary and sufficient
conditions under which the spectrum of the magnetic Schr\"odinger
operators  has an a.c. component.
\end{remark}

\no \textbf{Proof.} The determinant of $\big(\l I_\n-H_\a(k)\big)$
has the decomposition
\[\lb{Det}
\det\big(\l I_\n-H_\a(k)\big)=\prod_{j=1}^\n(\l-\l_{\a,j}(k))=
\l^\n+\x_1\l^{\n-1}+\x_2\l^{\n-2}+\ldots+\x_{\n-1}\l+\x_\n,
\]
where the coefficients $\x_j$ are given by (see, e.g., p. 87--88 in \cite{Ga60})
\[\lb{coDe}
\x_n=-{1\/n}\rt(T_n+\sum_{j=1}^{n-1}T_{n-j}\,\x_j \rt), \qqq T_n=\Tr
H^n_\a(k),\qqq n\in \N_\n,
\]
and, in particular, $\x_1=-T_1$, $\x_2=-{1\/2}(T_2-T_1^2),\,\ldots\,$.

(\emph{i}) $\Leftrightarrow$ (\emph{ii}). Let the spectrum of $H_\a$
be flat. Then all band functions $\l_{\a,j}(\cdot)$, $j\in\N_\n$,
are constant, and, consequently, $\Tr H_\a^n(k)=\sum_{j=1}^\n\l_{\a,j}^n(k)$ does not depend on $k$ for each $n\in\N_\n$.

Conversely, let $\Tr H_\a^n(k)$ does not depend on $k$ for each
$n\in\N_\n$. From this  and \er{Det}, \er{coDe} it follows that the
determinant $\det\big(\l I_\n-H_\a(k)\big)$ does not depend on $k$.
Then all band functions $\l_{\a,j}(\cdot)$, $j\in\N_\n$, are
constant and all spectral bands $\s_j(H_\a)$ are degenerate.

(\emph{ii}) $\Rightarrow$ (\emph{iii}). Let $\Tr H_\a^n(k)$ do not
depend  on $k$ for each $n\in\N_\n$. Then, using the Fourier series
\er{FsTrD}, we obtain
$$
\gh_{\a,n,\gm}=\sum_{\bc\in\wt\cC_{n,\gm}}\o(\bc)e^{-i\a(\bc)}=0,\qqq
\forall \,(n,\gm)\in\N_\n\ts\big(\Z^d\sm\{0\}\big).
$$

(\emph{iii}) $\Rightarrow$ (\emph{ii}). Let the condition \er{mpco}
hold true.  Then, by Theorem \ref{TFNL0}, we obtain
\[\lb{Trac}
\Tr H_\a^n(k)=\gh_{\a,n,0}=\sum_{\bc\in\wt\cC_{n,0}}\o(\bc)e^{-i\a(\bc)},
\qqq \forall\,n\in\N_\n,
\]
i.e., the traces $\Tr H_\a^n(k)$, $n\in\N_\n$, do not depend on $k$.

(\emph{iii}) $\Leftrightarrow$ (\emph{iv}). Using the Parseval's
identity for  the Fourier series \er{FsTrD}, we obtain
$$
\frac1{(2\pi)^d}\int_{\T^d}\Tr^2 H_\a^n(k)dk=
\sum_{\gm\in\Z^d}|\gh_{\a,n,\gm}|^2, \qqq \forall\, n\in\N.
$$
Thus, the condition \er{mpco} and \er{Trsq} are equivalent. \qq \BBox

\medskip

We prove Theorem \ref{Tfmp} about the flat
spectrum of the magnetic Schr\"odinger operators.

\medskip

\no \textbf{Proof of Theorem \ref{Tfmp}.} The first spectral band
$\s_1(H_0)$  of the Schr\"odinger operator $H_0$ without magnetic
field is non-degenerate (see Theorem 2.1.\emph{ii} in \cite{KS19}).
Then the a.c. spectrum of $H_0$ is not empty.
Then, by Theorem \ref{Tfbs}, there exists $(n,\gm)\in\N_\n\ts\big(\Z^d\sm\{0\}\big)$ such that the Fourier coefficient
$\gh_{0,n,\gm}\neq0$, where $\gh_{0,n,\gm}$ is defined by \er{FsTrD}. For this $(n,\gm)$ and for the magnetic potential $t\a$ we define the function
$$
f(t):=\gh_{t\a,n,\gm}=\sum_{\bc\in\wt\cC_{n,\gm}}\o(\bc)e^{-it\a(\bc)}, \qq t\in\R,
$$
and note that the sum is finite. This function has an analytic extension to the whole complex plane. If $f=\const$, then we obtain
$f=f(0)=\gh_{0,n,\gm}\ne0$ for such specific $(n,\gm)$. Then Theorem
\ref{Tfbs} yields that the a.c. spectrum of the operator $H_{t\a}$ is not empty for all $t\in \R$.

If $f\ne \const$, then $f$ has a finite number of zeros on any bounded
interval. Then Theorem \ref{Tfbs} yields that the a.c. spectrum of  $H_{t\a}$ is not empty for all except finitely many $t\in[0,1]$.
\qq \BBox

\smallskip

\no\textbf{Acknowledgments.} \footnotesize Our study was supported
by the   RFBR grant No. 19-01-00094.

%\medskip


\begin{thebibliography}{9999}
\setlength{\itemsep}{-\parskip}\footnotesize

\bibitem{BS98} Birman, M.Sh.; Suslina, T.A.  The two-dimensional periodic
 magnetic Hamiltonian is absolutely continuous, St. Petersburg Math. J., 9 (1998), no.~1, 21--32.

\bibitem{D85}  Danilov, L.I. On the spectrum of the Dirac
operator with periodic potential in $\R^n$, Theoret. and Math.
Phys., 85 (1990), no.1, 1039--1048.

\bibitem{Ga60} Gantmacher, F.R. The Theory of Matrices. NY: Chelsea Publishing, 1960.

\bibitem{HH95} Hempel, R.; Herbst, I. Bands and gaps for
periodic magnetic hamiltonians, Operator Theory: Advances and
Applications 78, Birkhauser, Basel, 1995, 175--184.

\bibitem{HS99} Higuchi, Y.; Shirai, T. The spectrum
of magnetic Schr\"odinger operators on a graph with periodic
structure, J. Funct. Anal., 169 (1999), 456--480.

\bibitem{HN09} Higuchi, Y.; Nomura, Y.
Spectral structure of the Laplacian on a covering graph.
European J. Combin. 30 (2009), no. 2, 570--585.

\bibitem{K97} Korotyaev, E. The estimates of periodic potentials in
terms of effective masses. Comm. Math. Phys. 183 (1997), no. 2,
383--400.

\bibitem{K00}  Korotyaev, E. Estimates for the Hill operator. I. J.
Differential Equations 162 (2000), no. 1, 1--26.

\bibitem{KS14} Korotyaev, E.; Saburova, N.
Schr\"odinger operators on periodic discrete graphs, J. Math. Anal. Appl.,
 420 (2014), no. 1, 576--611.

\bibitem{KS17}  Korotyaev, E.; Saburova, N.
Magnetic Schr\"odinger operators on periodic discrete graphs, J.
Funct. Anal., 272 (2017), 1625--1660.

\bibitem{KS19} Korotyaev, E.; Saburova, N. Spectral estimates for
 Schr\"odinger operators on periodic discrete graphs,
 St. Petersburg Math. J. 30 (2019), no.~4, 667--698.


\bibitem{RS78} Reed, M.; Simon, B. Methods of modern mathematical
physics, vol.IV. Analysis of operators, Academic Press, New York, 1978.

\bibitem{S99} Sobolev, A. Absolute continuity of the periodic magnetic
 Schr\"odinger operator. Invent. Math. 137 (1999), no. 1, 85--112.

\bibitem{T73} Thomas, L. Time dependent approach to scattering
 from impurities in a crystal, Comm. Math. Phys. 33 (1973), 335--343.

\end{thebibliography}
\end{document}